\documentclass{crm-article}

\usepackage{algorithm}
\usepackage{algorithmic}
\usepackage{float}
\usepackage{multirow}
\usepackage{hyperref}
 
\usepackage{ulem}

\begin{document}

\englishpaper 




\journalVol{10}
\journalNo{1} 
\setcounter{page}{1}

\journalSection{Математические основы и численные методы моделирования}
\journalSectionEn{Mathematical modeling and numerical simulation}

\journalReceived{00.00.2022.}
\journalAccepted{00.00.2022.}

\UDC{519.8}
\title{Подход к решению невыпуклой равномерно вогнутой седловой задачи со структурой}
\titleeng{An Approach for Non-Convex Uniformly Concave Structured Saddle Point Problem}
\thanks{Исследование выполнено при поддержке Министерства науки и высшего образования Российской Федерации (госзадание) № 075-00337-20-03, номер проекта 0714-2020-0005.}
\thankseng{The research is supported by the Ministry of Science and Higher Education of the Russian Federation (Goszadaniye) № 075-00337-20-03, project No. 0714-2020-0005.}

\author[1,2]{\firstname{М.\,С.}~\surname{Алкуса}}
\authorfull{Мохаммад С. Алкуса}
\authoreng{\firstname{M.\,S.}~\surname{Alkousa}}
\authorfulleng{Mohammad S. Alkousa}
\email{mohammad.alkousa@phystech.edu}
\affiliation[1]{Московский физико-технический институт, Россия, 141701, Московская~обл., г. Долгопрудный, Институтский пер.,~9}
\affiliationeng[1]{Moscow Institute of Physics and Technology, Russia, 141701, Moscow region, Dolgoprudny, Institutskiy per.,~9}

\affiliation[2]{Национальный исследовательский университет «Высшая школа экономики», Россия, 101000, г. Москва, ул. Мясницкая, д.~20}
\affiliationeng[2]{HSE University, Russia, 101000, Moscow, Myasnitskaya St.,~20}

\author[1,3,4]{\firstname{А.\,В. }~\surname{Гасников}}
\authorfull{Александр В. Гасников}
\authoreng{\firstname{A.\,V.}~\surname{Gasnikov}}
\authorfulleng{Alexander A. Gasnikov}
\email{gasnikov@yandex.ru}
\affiliation[3]{212705, г. Москва, Большой Каретный переулок, д.19 стр. 1, Институт проблем управления РАН им. А.А. Харкевича}
\affiliationeng[3]{Institute for Information Transmission Problems of the Russian Academy of Sciences (Kharkevich Institute)}
\affiliation[4]{3385000, Республика Адыгея, г. Майкоп, ул. Первомайск, д. 208, Кавказский математический центр Адыгейского государственного университета}
\affiliationeng[4]{Caucasus Mathematical Center, Adyghe State University, Maikop, Russia}

\author[5,3]{\firstname{П.\,Е.}~\surname{Двуреченский}}
\authorfull{Павел Е. Двуреченский }
\authoreng{\firstname{P.\,E.}~\surname{Dvurechensky}}
\authorfulleng{Pavel E. Dvurechensky}
\email{Pavel.Dvurechensky@wias-berlin.de}
\affiliation[5]{Weierstrass Institute for Applied Analysis and Stochastics, Берлин, 10117 Германия}
\affiliationeng[5]{Weierstrass Institute for Applied Analysis and Stochastics, Berlin, 10117 Germany}

\author[1]{\firstname{А.\,А.}~\surname{Садиев}}
\authorfull{Абдурахмон А. Садиев}
\authoreng{\firstname{A.\,A.}~\surname{Sadiev}}
\authorfulleng{Abdurakhmon A. Sadiev}
\email{sadiev.aa@phystech.edu}
  
\author[6]{\firstname{Л.\,Я.}~\surname{Разук}}
\authorfull{Лама Я. Разук}
\authoreng{\firstname{L.\,Y.}~\surname{Razouk}}
\authorfulleng{Lama Y. Razouk}
\email{lamarazouk94@gmail.com}
\affiliation[6]{Tartous University, Department of Mathematics, Tartous, Сирия}
\affiliationeng[6]{Tartous University, Department of Mathematics, Tartous, Syria}


\begin{abstract}
В последнее время седловым задачам уделяется большое внимание благодаря их мощным возможностям моделирования для множества задач из различных областей. Приложения этих задач встречаются во многочисленных современных прикладных областях, таких как робастная оптимизация, распределенная оптимизация, теория игр и приложения машинного обучения такие, как например минимизация эмпирического риска или обучение генеративно-состязательных сетей. Поэтому многие исследователи активно работают над разработкой численных методов для решения седловых задач в самых разных предположениях. Данная статья посвящена разработке численного метода решения седловых задач в невыпуклой равномерно-вогнутой постановке. В этой постановке считается, что по группе прямых переменных целевая функция может быть невыпуклой, а по группе двойственных переменных задача является равномерно вогнутой (это понятие обобщает понятие сильной вогнутости). Был изучен более общий класс седловых задач со сложной композитной структурой и Гёдьдерово непрерывными производными высшего порядка. Для решения рассматриваемой задачи был предложен подход, при котором мы сводим задачу к комбинации двух вспомогательных оптимизационных задач отдельно для каждой группы переменных: внешней задачи минимизации и внутренней задачи максимизации. Для решения внешней задачи минимизации мы используем \textit{адаптивный градиентный метод}, который применим для невыпуклых задач, а также работает с неточным оракулом, который генерируется путем неточного решения внутренней задачи максимизации. Для решения внутренней задачи максимизации мы используем \textit{обобщенный ускоренный метод с рестартами}, который представляет собой метод, объединяющий методы ускорения высокого порядка для минимизации выпуклой функции, имеющей Гёльдерово непрерывные производные высшего порядка. Важной компонентой проведенного анализа сложности предлагаемого алгоритма является разделение оракульных сложностей на число вызовов оракула первого порядка для внешней задачи минимизации и оракула более высокого порядка для внутренней задачи максимизации. Более того,  оценивается сложность всего предлагаемого подхода.
\end{abstract}

\begin{abstracteng}
Recently, saddle point problems have received much attention due to their powerful modeling capability for a lot of problems from diverse domains. Applications of these problems occur in many applied areas, such as robust optimization, distributed optimization, game theory, and many applications in machine learning such as empirical risk minimization and generative adversarial networks training. Therefore, many researchers have actively worked on developing numerical methods for solving saddle point problems in many different settings. This paper is devoted to developing a numerical method for solving saddle point problems in the non-convex uniformly-concave setting. We study a general class of saddle point problems with composite structure and H\"older-continuous higher-order derivatives. To solve the problem under consideration, we propose an approach in which we reduce the problem to a combination of two auxiliary optimization problems separately for each group of variables, outer minimization problem w.r.t. primal variables, and inner maximization problem w.r.t the dual variables. For solving the outer minimization problem, we use the  \textit{Adaptive Gradient Method}, which is applicable for non-convex problems and also works with an inexact oracle that is generated by approximately solving the inner problem. For solving the inner maximization problem, we use the \textit{Restarted Unified Acceleration Framework}, which is a framework that unifies the high-order acceleration methods for minimizing a  convex function that has  H\"older-continuous higher-order derivatives. Separate complexity bounds are provided for the number of calls to the first-order oracles for the outer minimization problem and higher-order oracles for the inner maximization problem. Moreover, the complexity of the whole proposed approach is then estimated.
\end{abstracteng}
\keywordeng{Saddle point problem}
\keywordeng{non-convex optimization}
\keywordeng{uniformly convex function}
\keywordeng{inexact oracle}
\keywordeng{higher-order method}

\maketitle

\paragraph{1. Introduction}\label{sec_introduction}
Due to their numerous applications, saddle point (or minimax) problems have recently gained significant interest in many fields, including  robust optimization \cite{ben2009robust,namkoong2016stochastic}, empirical risk minimization \cite{tan2018stochastic,zhang2015stochastic},  statistics \cite{berger2013statistical}, distributed non-convex optimization \cite{lu2019block}, game theory \cite{myerson2013game,nouiehed2019solving} and many applications of machine learning, including reinforcement learning \cite{dai2018sbeed}, adversarial learning \cite{sinha2017certifying}, learning exponential families \cite{dai2019kernel}, fair statistical inference \cite{madras2018learning,xu2018fairgan}, generative adversarial networks \cite{zhang2021understanding,lei2020sgd,sanjabi2018}, domain adversarial training \cite{zhang2021understanding}, Wasserstein robust models \cite{sinha2017certifying}, robust learning over multiple domains \cite{qian2019robust}, off-policy reinforcement learning \cite{huang2021convergence}, prediction and regression problems \cite{taskar2005structured}.
Most of these applications require solving a saddle point problem in which the objective function is non-convex in one variable and strongly-concave (or in general uniformly-concave) in another variable. For the saddle point problems, there is a long history of studies. Recently, in the last years, many researchers have actively worked on the development of numerical methods for solving these problems in many different settings.

The general setting of saddle point problems, which is Non-Convex Concave (without strong concavity assumption), is extensively explored in a lot of works \cite{kong2019accelerated,nouiehed2019solving,ostrovskii2020efficient,thekumparampil2019efficient,zhang2020single}. The work \cite{nouiehed2019solving},  was the first providing non-asymptotic convergence rates for non-convex concave problems without assuming special structure of the objective function. The authors of \cite{thekumparampil2019efficient} proposed an algorithm for the general smooth non-convex concave saddle point problem with complexity $\widetilde{\mathcal{O}} \left(\varepsilon^{-3}\right)$.\footnote{Here and below we use the notation $\mathcal{O}(\cdot)$ to denote non-asymptotic upper bounds up to constant factors. More precisely, $a=\mathcal{O}(b)$ means that there exist a constant $C$ such that $a \leqslant C b$. We use $\widetilde{\mathcal{O}}(\cdot)$ for the same purposes, but when the constant $C$ may include polylogarithmic factors of the parameters of interest, e.g., the desired accuracy $\varepsilon$ and confidence probability $\sigma$.}  The same complexity for a simple two-time-scale inexact proximal point algorithm was proved in \cite{yang2020catalyst}.  In  \cite{zhang2020single}, a single-loop algorithm ''\textit{Smoothed GDA}'' (which is an iterative algorithm where each iteration step has a closed form update) was proposed with iteration complexity $\mathcal{O}\left(\varepsilon^{-2}\right)$ for the special case of the general non-convex concave problems, which is equivalent to the problem of minimizing the point-wise maximum of a finite collection of functions. The complexity  $\mathcal{O}\left(\varepsilon^{-2}\right)$ improves other complexities, such as $\mathcal{O}\left(\varepsilon^{-2.5}\right)$ for triple-loop algorithms \cite{ostrovskii2020efficient,lin2020nearoptimal} and $\mathcal{O}\left(\varepsilon^{-3.5}\right)$ for double-loop algorithms \cite{nouiehed2019solving}. Also, non-convex concave saddle point problems via zeroth-order algorithms  were studied in \cite{xu2021zeroth}, where the authors proposed a zeroth-order alternating randomized gradient projection algorithm for smooth non-convex concave problems. 

Under the strong concavity assumption,  several algorithms were proposed for the non-convex-strongly-concave  (NC-SC) class of problems, including GDmax \cite{nouiehed2019solving}, GDA \cite{lin2020gradient}, alternating GDA \cite{boct2020alternating,xu2020unified,yang2020global}. Recently, in \cite{lin2020nearoptimal}, the authors proposed an accelerated algorithm for smooth NC-SC and  non-convex concave minimax problems. This algorithm achieves the best dependency on the condition number by combining proximal point algorithm with accelerated gradient descent. Namely, it achieves a gradient complexity bound of $\widetilde{\mathcal{O}}\left(\sqrt{\kappa} \varepsilon^{-2} \right)$ which improves upon the best known bound $\widetilde{\mathcal{O}}\left(\kappa^2 \varepsilon^{-2}\right)$ \cite{lin2020gradient}. Another line of research is devoted to solving  NC-SC minimax problems under additional structural properties \cite{diakonikolas2021efficient,lin2018solving,song2021optimistic,yang2020global}. For more details in the recent works in NC-SC in different settings, including the general and average-smooth finite-sum settings, see \cite{zhang2021complexity,luo2020stochastic,ostrovskii2021nonconvex,luo2021finding} and references therein.

Despite this active line of research for different settings of saddle point problems, to the best of our knowledge, the case when the objective function is non-convex in one variable and uniformly concave in another variable instead of strong concavity assumption (where the strong concavity case is a spacial case of the inform concavity, see Definition \ref{def_uniformly_convex} below),  still are not investigated. Therefore, we  address this setting of the problem in this paper that is devoted to the development of a numerical method for solving saddle point problems in the non-convex uniformly-concave setting. Moreover, we study a more general class of saddle point problems with composite structure and H\"older-continuous higher-order derivatives. In order to solve the problem under consideration, we propose an approach, in which we reduce the problem to a combination of two auxiliary optimization problems separately for each group of variables: outer and inner problems. For solving the outer minimization problem, we use the  \textit{Adaptive Gradient Method} \cite{Dvurechensky_2017}, which is applicable for non-convex problems and also works with inexact oracle. For solving the inner maximization problem, we use the \textit{Restarted Unified Acceleration Framework} \cite{song2019unified}, which is a framework, that unifies the high-order  acceleration  methods  for  minimizing  a  convex  function  that  has  H\"older-continuous higher-order derivatives. Separate bounds are provided for the number of calls to the first-order oracles for the outer minimization problem and higher-order oracles for inner maximization problem. Moreover, the complexity of the whole proposed approach is then estimated.

The paper consists of the introduction, conclusion, and two main sections.  In Section 2, we present the problem statement, its connected assumptions, and some necessary facts that will be used in the proposed approach. In Section 3, we present the proposed approach in order to solve the considered saddle point problem by applying \textit{Adaptive Gradient Method} \cite{Dvurechensky_2017} and \textit{Restarted Unified Acceleration Framework} \cite{song2019unified}. When analyzing the proposed approach, separate complexity bounds are provided for the number of calls to the first-order oracles for the outer minimization and higher-order oracles for the inner maximization problems. Moreover, the complexity of the whole proposed approach is then estimated.

\paragraph{2. Preliminaries and Problem Statement}\label{section_preliminaries}

\subparagraph{Problem Formulation and Assumptions}
We focus on the following structured composite saddle point problem
\begin{gather}\label{main_spp_problem}
    \min_{x \in \mathcal{X}}\max_{y \in \mathcal{Y}}\left\{S(x,y) := F(x,y) -h(y) +r(x)\right\},
\end{gather}
in the case when $\mathcal{X} \subseteq \mathbb{R}^{n_x}$ is a compact convex set, and $\mathcal{Y} =  \mathbb{R}^{n_y}$. The function $h: \mathcal{Y} \to \mathbb{R}$ is  differentiable and uniformly convex (see Definition \ref{def_uniformly_convex}, below), the function $r: \mathcal{X} \to \mathbb{R}$ is a simple (maybe non-smooth) convex function, e.g. $\|x\|_1$. Note that when $h(y) = r(x) = 0, \; \forall x \in \mathcal{X}, \forall y \in \mathcal{Y}$, we obtain a problem which is a special (unstructured) case of problem \eqref{main_spp_problem}. 
We assume that the function   $F(\cdot,y)$, for every $y \in \mathcal{Y}$, is possibly non-convex in the primal variable $x \in \mathcal{X}$, and $F(x,\cdot)$ for every $x \in \mathcal{X}$ is concave in the dual variable $y \in \mathcal{Y}$. Thus the function $S(x,y)$ will be non-convex uniformly concave (NC-UC), for every $x \in \mathcal{X}$ and $y \in \mathcal{Y}$. Moreover, we assume that the $p$-th derivatives of $F$ with respect to its variables  satisfy the H\"{o}lder condition (in other words, the function $F$ has $(p, \nu)$--H\"{o}lder-continuous derivatives, for some $\nu \in [0,1], p \in \{1,2, \ldots\}$).  More precisely, for arbitrary $x,x' \in \mathcal{X}$, $y,y' \in \mathcal{Y}$, the following inequalities hold (here and in what follows, the norm $\|\cdot\|$ denotes  the Euclidean  norm)
\begin{gather}\label{smooth_F_1}
  \left\|\nabla_x^p F(x,y) - \nabla_x^p F(x', y)\right\| \leqslant L_{xx}\|x -x'\|^{\nu},
\end{gather}
\begin{gather}\label{smooth_F_2}
    \left\|\nabla_x^p F(x,y) - \nabla_x^p F(x, y')\right\| \leqslant L_{xy}\|y -y'\|^{\nu},
\end{gather}
\begin{gather}\label{smooth_F_3}
   \left\|\nabla_y^p F(x,y) - \nabla_y^p F(x', y)\right\| \leqslant L_{yx}\|x -x'\|^{\nu},
\end{gather}
\begin{gather}\label{smooth_F_4}
    \left\|\nabla_y^p F(x,y) - \nabla_y^p F(x, y')\right\| \leqslant L_{yy}\|y -y'\|^{\nu},
\end{gather}
where $L_{xx}, L_{yy}, L_{xy},  L_{yx}> 0$.

Let
\begin{gather}\label{function:S_hat}
    \widehat{S}(x,y) := F(x,y) - h(y).
\end{gather}
and
\begin{gather}\label{inner_max_problem}
    g(x) := \max_{y \in \mathcal{Y}} \widehat{S}(x,y).
\end{gather}
Then, the problem \eqref{main_spp_problem} can be rewritten as follows
\begin{gather}\label{equivalent_min_problem}
    \min_{x \in \mathcal{X}} \{g(x) +r(x)\}.
\end{gather}

Let 
\begin{gather}\label{solution:max_problem_g(x)}
    y^*(x) \in \arg\max_{y \in\mathcal{Y}} \widehat{S}(x,y) \quad
     \forall x \in \mathcal{X}
\end{gather}
be a solution of the maximization problem \eqref{inner_max_problem}. Then 
\begin{gather}\label{function:g(x)}
    g(x) = \widehat{S}\left(x,y^*(x)\right) = F\left(x,y^*(x)\right) - h\left(y^*(x)\right).
\end{gather}

 For an arbitrary $x \in \mathcal{X}$ and some $\hat{\delta} \geqslant 0$, the point $\widetilde{y}_{\hat{\delta}}(x) \in \mathcal{Y}$ is called a $\hat{\delta}$-approximate solution of problem \eqref{inner_max_problem} if
\begin{gather*}\label{delta_inexact_solution}
g(x) - \widehat{S}\left(x,\widetilde{y}_{\hat{\delta}}(x) \right) = \widehat{S}(x,y^*(x) ) - \widehat{S}\left(x,\widetilde{y}_{\hat{\delta}}(x) \right) \leqslant \hat{\delta}.
\end{gather*}

\subparagraph{Necessary Auxiliary Statements }
In order to solve the considered problem \eqref{main_spp_problem}, we propose an approach that uses several auxiliary results and algorithms as building blocks. In what follows, we present these auxiliary results, whereas the used algorithms in the proposed approach will be presented in the next section.

\begin{fed}\label{def_uniformly_convex}
We say that differentiable function $f$ is uniformly convex of degree $q \geqslant 2$ on a convex set $Q \subseteq dom f$  if for some constant $\sigma_q >0$ it satisfies the following inequality:
\begin{gather*}\label{eq_uniformly_convex}
    f(y) \geqslant f(x)+\langle\nabla f(x), y-x\rangle+\frac{\sigma_{q}}{q}\|y-x\|^{q}, \quad  \forall x, y \in Q.
\end{gather*}
We note that the uniformly convex functions of degree $q = 2$ are strongly convex with $\sigma_2$ as the strong convexity parameter.
\end{fed}

The next lemma gives us the main properties of the function $g$ defined in \eqref{inner_max_problem} and the mapping  $y^*$  defined in \eqref{solution:max_problem_g(x)}.

\begin{lem}\label{lemma_function_g}
Let $h:\mathbb{R}^{n_y} \to \mathbb{R}$ be a uniformly convex function of degree $q \geqslant 2$ with constant $\sigma_q>0$, and the function $F$ have $(1, \nu)$-H\"older-continuous derivative with respect to the variable $x \in \mathcal{X}$, for $\nu \in [0,1]$. Then the mapping $y^*(\cdot)$ satisfies the H\"{o}lder condition on the set $\mathcal{X}$. Also, $\nabla g$ is H\"{o}lder-continuous with H\"{o}lder constant $L_{\nu}(g) := L_{xy} \left(\frac{qL_{xy}}{\sigma_q}\right)^{\frac{\nu}{q-\nu}} + L_{xx} D_0^{\frac{\nu(q- \nu-1)}{q-\nu}}$ and H\"{o}lder exponent $\nu_g := \frac{\nu}{q-\nu}$, where $D_0= \sup_{x,x'\in\mathcal{X}}\|x' -x\|$. 
\end{lem}
\textbf{Proof:}
Similarly to \cite{Alkousa_2020}. Since $\widehat{S}(x_1, \cdot)$ is  uniformly concave of degree $q$ with constant $\sigma_q$, for arbitrary $x_1, x_2 \in \mathcal{X}$, we have
\begin{gather*}\label{ineq:1}
   \|y^*(x_1) - y^*(x_2)\|^q \leqslant \frac{q}{\sigma_q} \left( \widehat{S}(x_1, y^*(x_1)) - \widehat{S}(x_1, y^*(x_2)) \right).
\end{gather*}

On the other hand, since  $y^*(x_2)$ attains the maximum to $\widehat{S}(x_2,\cdot )$, then $\widehat{S}(x_2, y^*(x_1)) - \widehat{S}(x_2, y^*(x_2)) \leqslant 0$. Therefore, we can write 
\begin{gather*}\label{ineq:2}
\begin{split}
\widehat{S}(x_1, y^*(x_1)) &-  \widehat{S}(x_1, y^*(x_2))  \leqslant  \widehat{S}(x_1, y^*(x_1)) - \widehat{S}(x_1, y^*(x_2))  -  \widehat{S}(x_2, y^*(x_1)) + \widehat{S}(x_2, y^*(x_2)) = \\
& = \left( F(x_1, y^*(x_1)) - F(x_1, y^*(x_2)) \right) - \left( F(x_2, y^*(x_1)) - F(x_2, y^*(x_2)) \right) \\
& = \int_0^1 \langle  \nabla_xF(x_1 + t(x_2 - x_1), y^*(x_1)) - \nabla_xF(x_1 + t(x_2 - x_1), y^*(x_2)), x_2 - x_1 \rangle dt  \\
&  \leqslant  \| \nabla_xF(x_1 + t(x_2 - x_1), y^*(x_1)) - \nabla_xF(x_1 + t(x_2 - x_1), y^*(x_2))  \| \cdot \|x_2 - x_1\|  \\
&   \leqslant L_{xy} \|y^*(x_1) - y^*(x_2)\|^{\nu} \cdot \|x_2 - x_1\|,
\end{split}
\end{gather*}
where in the last inequality we used \eqref{smooth_F_2} for $p = 1$. Thus, we have the following inequality
\begin{gather*}\label{ineq:Lip_y_star}
    \|y^*(x_2) - y^*(x_1)\| \leqslant \left(\frac{q L_{xy}}{\sigma_q}\right)^{\frac{1}{q-\nu}}\|x_2 - x_1\|^{\frac{1}{q-\nu}}, \quad \forall x_1, x_2 \in \mathcal{X}, 
\end{gather*}
which means that $y^*(\cdot)$ satisfies the H\"{o}lder condition on $\mathcal{X}$.

Similarly to \cite{nouiehed2019solving}, we find
\begin{gather}\label{nabla_g}
    \nabla g(x) = \nabla_x \widehat{S}(x,y^*(x)) = \nabla_x F(x,y^*(x)), \quad \forall x \in\mathcal{X}.
\end{gather}

Further, from \eqref{nabla_g}, we have
\begin{align*}
\|\nabla g(x_1) &- \nabla g(x_2)\| = \|\nabla_x F(x_1,y^*(x_1))-\nabla_x F(x_2,y^*(x_2))\| = \\
&  = \|\nabla_x F(x_1,y^*(x_1)) - \nabla_x F(x_1,y^*(x_2)) + \nabla_x F(x_1,y^*(x_2)) - \nabla_x F(x_2,y^*(x_2))\|  \\
&  \leqslant  \|\nabla_x F(x_1,y^*(x_1)) - \nabla_x F(x_1,y^*(x_2))\| + \|\nabla_x F(x_1,y^*(x_2)) - \nabla_x F(x_2,y^*(x_2))\|  \\
&  \leqslant L_{xy}\|y^*(x_1)- y^*(x_2)\|^{\nu} + L_{xx} \|x_1 - x_2\|^{\nu}  \\
& \leqslant L_{xy} \left(\frac{q L_{xy}}{\sigma_q}\right)^{\frac{\nu}{q-\nu}}\|x_1 - x_2\|^{\frac{\nu}{q-\nu}} + L_{xx}\|x_1 - x_2\|^{\nu} \\
& = L_{xy} \left(\frac{q L_{xy}}{\sigma_q}\right)^{\frac{\nu}{q-\nu}}\|x_1 - x_2\|^{\frac{\nu}{q-\nu}} + L_{xx}\|x_1 - x_2\|^{\frac{\nu}{q-\nu}}\cdot \|x_1 - x_2\|^{\frac{\nu(q-\nu -1)}{q-\nu}}  \\
& \leqslant \left(L_{xy}\left(\frac{qL_{xy}}{\sigma_q}\right)^{\frac{\nu}{q-\nu}} + L_{xx} D_0^{\frac{\nu(q- \nu-1)}{q-\nu}}\right)\cdot\|x_1-x_2\|^{\frac{\nu}{q-\nu}}.
\end{align*}

This means, that  $\nabla g$ is H\"{o}lder-continuous with the H\"{o}lder constant
$$
   L_{\nu}(g) = L_{xy} \left( \frac{qL_{xy}}{\sigma_q}\right)^{\frac{\nu}{q-\nu}} + L_{xx} D_0^{\frac{\nu(q- \nu-1)}{q-\nu}},
$$ 
and the H\"{o}lder exponent $\nu_g = \frac{\nu}{q-\nu}$.
\qed

\paragraph{3. Used Algorithms and The Proposed Approach}\label{fund_section}
When solving the considered saddle point problem \eqref{main_spp_problem}, as we saw, we deal with two problems, the outer is a minimization problem and the inner is a  maximization problem. Using an iterative method for the outer problem requires solving the inner problem numerically in each iteration. Let us first mention some methods, that will be used in order to solve the outer and inner optimization problems.

\subparagraph{Algorithm for the outer minimization problem:}
For the outer minimization problem \eqref{equivalent_min_problem}, we will use an algorithm that was proposed in \cite{Dvurechensky_2017,bogolubsky2016learning}. This algorithm   (listed as Algorithm \ref{method_Pavel_2017} below) is a first-order method and it was developed for composite non-convex minimization problems with inexact oracle. Moreover, it is universal with respect to H\"older parameters of the problem. Firstly, let us introduce some fundamental concepts, connected with Algorithm \ref{method_Pavel_2017}.

Let $d : \mathcal{X} \rightarrow \mathbb{R}$ be a prox-function (distance generating function), i.e. it is continuous, convex
on $\mathcal{X}$ and
\begin{enumerate}
    \item admits a continuous selection of subgradients $\nabla d(x)$, in $x \in \mathcal{X}^{\circ}$,  where $ \mathcal{X}^{\circ} \subseteq \mathcal{X}$ is the set of all $x$ such that $\nabla d(x)$ exists.
    \item $d$ is $1$-strongly convex with respect to the norm $\|\cdot\|$, i.e.
\begin{gather*}
	d(y) \geqslant d(x) +\langle \nabla d(x), y-x \rangle + \frac{1}{2} \| y-x \|^2, \quad \forall\;  x\in \mathcal{X}^{\circ}, y \in \mathcal{X}.
\end{gather*}
\end{enumerate}

The corresponding Bregman divergence is defined as
$$
V_z(x) = d(x) - d(z) - \left\langle \nabla d(z), x-z\right\rangle, \quad \forall x \in \mathcal{X}, z \in \mathcal{X}^{\circ}. 
$$

In particular, in the standard proximal setup (i.e. Euclidean setup) we can choose $d(x) = \frac{1}{2}\|x\|_2^2$ (Euclidean prox function), leading to $V_z(x) = \frac{1}{2}\|x - z\|_2^2$.

Now, let we consider the following  minimization problem
\begin{gather}\label{non_convex_comp_problem}
    \min_{x \in \mathcal{X}}\{f(x) = \xi(x) + \zeta(x)\},
\end{gather}
where $\zeta(x)$ is a simple convex function, and $\xi(x)$ is a non-convex function, endowed with an inexact first-order oracle (see Definition \ref{def_ineexact_oracle}, below), and has H\"{o}lder-continuous gradient on $\mathcal{X}$, with constant $L_{\nu}$ and exponent $\nu$.

\begin{fed}\label{def_ineexact_oracle}(\cite{Dvurechensky_2017}) 
We say that a function $\xi(x)$ is equipped with an inexact first-order oracle on a set $\mathcal{X}$ if there exists $\delta_u > 0$ and at any point $x \in \mathcal{X}$ for any number $\delta_c > 0$ there exists a constant $L(\delta_c) \in (0, +\infty)$ and one can calculate $\widetilde{\xi}(x, \delta_c, \delta_u)$ and $\widetilde{g}(x, \delta_c, \delta_u)$ satisfying 
$$
\left|\xi(x) - \widetilde{\xi}(x, \delta_c, \delta_u)\right| \leqslant \delta_c + \delta_u,
$$
and
$$
\xi(y) - \left(\widetilde{\xi}(x, \delta_c, \delta_u) - \langle \widetilde{g}(x, \delta_c, \delta_u), y - x\rangle\right) \leqslant \frac{L(\delta_c)}{2} \|x - y\|^2 + \delta_c + \delta_u, \quad \forall y \in \mathcal{X}.
$$

\end{fed}

\begin{rem} 
By Lemma \ref{lemma_function_g}, we have
$$
\|\nabla g(x_1) - \nabla g(x_2)\| \leqslant L_{\nu}(g) \|x_1 - x_2 \|^{\nu_g}, \quad \forall x_1, x_2 \in \mathcal{X}.
$$
Then 
$$
g(z) \leqslant g(x) + \langle \nabla g(x), z - x \rangle + \frac{L_{\nu}(g)}{1+ \nu_g} \|x - z\|^{1+ \nu_g}, \quad \forall x, z \in \mathcal{X}.
$$
Now, for all $x \in \mathcal{X}$ and $\delta > 0$, we have (see \cite{nesterov2015universal}, Lemma 2), 
$$
g(z) - \left( g(x) - \langle \nabla g(x), z - x\rangle \right) \leqslant \frac{L(\delta)}{2}\|x - z\|^2 + \delta, \quad \forall z \in \mathcal{X},
$$
where 
\begin{gather}\label{L_delta_g}
     L(\delta) = \left( \frac{1- \nu_g}{1+\nu_g}\cdot \frac{2}{\delta} \right)^{\frac{1- \nu_g}{1+\nu_g}} \cdot \left( L_{\nu}(g) \right)^{\frac{2}{1+\nu_g}}.
\end{gather}
Thus, according to the Definition \ref{def_ineexact_oracle}, $\left(g(x), \nabla g(x)\right)$ is an inexact first-order oracle with $\delta_u = 0, \delta_c  = \delta$, and $L(\delta)$ given by \eqref{L_delta_g}. Note that, if $\left(g(x), \nabla g(x)\right)$ can only be calculated inexactly, then their approximations will again be an inexact first-order oracle. 
\end{rem}

\begin{fed}\label{def_2_Pavel_paper}(\cite{Dvurechensky_2017})
Let $\mathbb{E}$ be a finite-dimensional real vector space and $\mathbb{E}^*$ be its dual. Assume that we are given $\delta_{pc} >0, \gamma>0,\bar{x} \in \mathcal{X}^{\circ}$ and $\eta \in \mathbb{E}^*$. We call a point $\widetilde{x} = \widetilde{x}(\bar{x}, \eta, \gamma, \delta_{pc}, \delta_{pu}) \in \mathcal{X}^{\circ}$, an inexact composite prox-mapping iff for any $\delta_{pc} >0$ we can calculate $\widetilde{x}$ and there exists $s \in \partial \zeta(\widetilde{x})$ such that it holds that
\begin{gather*}
   \left\langle \eta + \frac{1}{\gamma} \left( \nabla d(\widetilde{x}) - \nabla d(\bar{x})\right) + s, z - \widetilde{x}  \right\rangle \geqslant -\delta_{pc}-\delta_{pu}, \quad \forall z \in \mathcal{X},
\end{gather*}
and we write
\begin{gather*}
    \widetilde{x} = {\arg\min_{x \in \mathcal{X}}}^{\delta_{pc}+\delta_{pu}} \left\{\langle \eta,x\rangle+\frac{1}{\gamma}V_{\bar{x}}(x) +\zeta(x)\right\}.
\end{gather*}
\end{fed}

\begin{algorithm}[!ht]
\caption{Adaptive Gradient Method for Problems with Inexact Oracle \cite{Dvurechensky_2017}.}
\label{method_Pavel_2017}
	\begin{algorithmic}[1]
		\REQUIRE starting point $x_0 \in \mathcal{X}^{\circ}$, accuracy $\varepsilon >0$, initial guess $L_0 > 0$, $\delta_u>0$ and $\delta_{pu}>0$.
		\STATE Set $k=0$.
		\REPEAT
		\STATE Set $M_k = L_k/2$.
		\REPEAT
		\STATE Set $M_k = 2M_k, \delta_{c,k} = \delta_{pc,k} = \frac{\varepsilon}{20M_k}$.
		\STATE Calculate  $\widetilde{\xi}\left(x_{k}, \delta_{c, k}, \delta_{u}\right)$ and $\widetilde{g}\left(x_{k}, \delta_{c, k}, \delta_{u}\right)$.
		\STATE Calculate
		\begin{gather*}
            z_{k}={\arg \min_{x \in \mathcal{X}}}^{\delta_{p c, k}+\delta_{p u}}\left\{\left\langle\widetilde{g}\left(x_k, \delta_{c, k}, \delta_u\right), x\right\rangle + M_k V_{x_k}(x)+\zeta(x)\right\}.
        \end{gather*}
		\STATE Calculate $\widetilde{\xi}\left(z_k, \delta_{c, k}, \delta_u\right)$.
		\UNTIL{
		\begin{gather*}
	    \widetilde{\xi}\left(z_{k}, \delta_{c, k}, \delta_{u}\right) \leqslant \widetilde{\xi}\left(x_k, \delta_{c, k}, \delta_u\right)+\left\langle\widetilde{g}\left(x_k, \delta_{c, k}, \delta_u\right), z_k - x_k\right\rangle + \frac{M_k}{2}\left\|z_k - x_{k}\right\|^{2}+
		+\frac{\varepsilon}{10 M_k}+2 \delta_u.
        \end{gather*}
		}
		\STATE Set $x_{k+1} = z_k, L_{k+1} = M_k/2$ and $k=k+1$.
		\UNTIL{$\min_{i \in 1, \ldots, k}\left\|M_{i}\left(x_i - x_{i+1}\right)\right\| \leqslant \varepsilon$.
		}
		\ENSURE The point $x_{K+1}$, such that $K=\arg \min_{i \in 1, \ldots, k}\left\|M_i\left(x_i - x_{i+1}\right)\right\|$.
	\end{algorithmic}	
\end{algorithm}

In \cite{Dvurechensky_2017}, for Algorithm \ref{method_Pavel_2017}, it was proved that if $L(\delta_c)$ in Definition \ref{def_ineexact_oracle}, is given by
\begin{gather*}
L\left(\delta_c\right) = \left(\frac{1-\nu}{1+\nu} \cdot \frac{2}{\delta_c}\right)^{\frac{1-\nu}{1+\nu}} L_\nu^{\frac{2}{1+\nu}}, \quad \delta_c > 0, \;\text{and} \; \nu \in (0,1].
\end{gather*}
Then after 
\begin{gather}\label{complexity_alg_adaptive_grad}
    \mathcal{O}\left( \frac{L_{\nu}^{\frac{1}{\nu}}\left(f(x_0)-f^*\right)}{\varepsilon^{\frac{1+\nu}{2 \nu}}}\right)
\end{gather}
iterations of Algorithm \ref{method_Pavel_2017},  it holds that
\begin{gather*}
    \left\|M_k(x_K - x_{K+1})\right\|^2 \leqslant \varepsilon.
\end{gather*}

\subparagraph{Algorithm for the inner maximization problem:}

For the inner maximization problem \eqref{inner_max_problem}, when $\mathcal{Y} = \mathbb{R}^{n_y}$, we will use the Restarted Unified Acceleration Framework ({\tt Restarted UAF}) \cite{song2019unified}, see Algorithm \ref{restart_UAM} below.  The {\tt Restarted UAF} algorithm was proposed in \cite{song2019unified} for uniformly convex functions, it represents a general restart scheme for such general class of problems. The Unified Acceleration Framework ({\tt UAF}) \cite{song2019unified}, is a framework, that unifies the high-order acceleration methods for minimizing a convex function that has H\"older-continuous derivatives. The iteration complexities of instances of both the {\tt UAF} and the {\tt Restarted UAF} match existing lower bounds in most important cases \cite{grapiglia2019tensor}. 
Let us here mention briefly the statements of the algorithm {\tt Restarted UAF}, in order to solve the following composite optimization problem
\begin{gather}\label{composite_problm_meta}
    \min_{y \in \mathbb{R}^{n_y}} \{f(y) = \varphi(y)+\psi(y)\},
\end{gather}
where $\varphi, \psi$  are convex functions. Assume that
\begin{enumerate}
    \item \textbf{Assumption 1}: the function $\varphi$ has $(p, \nu)$-H\"older-continuous derivatives, with the constant of smoothness $L>0$. This means
    $$
    \left\|\nabla^{p} \varphi(y)-\nabla^{p} \varphi(y^{\prime})\right\| \leqslant L\|x-y\|^{\nu}, \quad \forall y, y^{\prime} \in \mathbb{R}^{n_y},
    $$
    where $\nu \in [0,1]$  and $p\in \{1,2,\ldots\}$.
    
    \item \textbf{Assumption 2}: the objective function $f$ is uniformly convex of degree $q \geqslant 2$ and constant $\sigma_q >0$.
\end{enumerate}

For the problem \eqref{composite_problm_meta}, let $\mathcal{A}_m(y)$ denote the output of an algorithm $\mathcal{A}$ after $m$ iterations with an input $y$, which satisfies
    \begin{gather}\label{1}
        f\left(\mathcal{A}_{m}(y)\right)-f\left(y_{*}\right) \leqslant \frac{c_{\mathcal{A}}\left\|y-y_{*}\right\|^{\nu}}{m^{r}},
    \end{gather}
    for some constants $r >0, \nu >0, c_{\mathcal{A}}>0$, and $y_*$ is a solution of \eqref{composite_problm_meta}.
    
Let $R >0$ be a constant such that $\|y_0 - y_*\|  \leqslant R$, where $y_0$ is the starting point of algorithm $\mathcal{A}$. We define
\begin{gather}\label{2}
m_{0}=\left[\left(\frac{2^{q} q c_{\mathcal{A}} R^{\nu -q}}{\sigma_q}\right)^{\frac{1}{r}}\right], \quad k_{0}= \begin{cases}{\left[\frac{1}{q}+\frac{\nu}{q} \log _{2} (R)+\frac{1}{\nu-q} \log_{2}\left(\frac{q c_{\mathcal{A}}}{\sigma_q}\right)\right],} &  q< \nu \\ +\infty \qquad\qquad\qquad\qquad\qquad\qquad\quad \, ,& q \geqslant  \nu \end{cases}
\end{gather}

\begin{algorithm} [H]
\caption{  Restarted Unified Acceleration Framework ({\tt Restarted UAF}) \cite{song2019unified}.
}
\label{restart_UAM}
	\begin{algorithmic}[1]
		\REQUIRE $q \geqslant 2$, starting point $y_0 \in \mathbb{R}^{n_y}$, $K \in \mathbb{Z}_+$, an algorithm $\mathcal{A}$ satisfying \eqref{1}, constants $m_0$ and $k_0$ which defined in \eqref{2}. 
		\STATE Set $z_0 = y_0$.
		\FOR{ $k = 0, 1, \ldots, K - 1 $}
		\IF{$k \leqslant k_0 - 1$}
		\STATE $m_{k}=\left[m_{0} 2^{-\frac{\nu-q}{r} k}\right]$.
		\STATE $z_{k+1}=\mathcal{A}_{m_{k}}\left(z_{k}\right)$.
		\ELSE
		\STATE  $z_{k+1}=\mathcal{A}_1\left(z_{k}\right)$. \ENDIF
		\ENDFOR
		\ENSURE $y_{K} := z_{K}$.
	\end{algorithmic}
\end{algorithm}

In \cite{song2019unified}, for the Algorithm \ref{restart_UAM}, it was proved the following result.

\begin{teo}\label{theorem_restart_UAM}
Under assumptions 1 and 2 above, to achieve an $\varepsilon$-solution of the problem \eqref{composite_problm_meta} by Algorithm \ref{restart_UAM} with $\mathcal{A}$ being {\tt UAF}, the number of iterations we need is at most 
\begin{gather*}\label{complexity1}
\mathcal{O}\left(\left(\frac{L}{\sigma_q}\right)^{\frac{2}{3(p+\nu)-2}} \log \left(\frac{1}{\varepsilon}\right)\right) \quad \text{if} \;\; q=p+\nu,
\end{gather*}
\begin{gather*}\label{complexity2}
\mathcal{O}\left(\left(\frac{L}{\sigma_q}\right)^{\frac{2}{3(p+\nu)-2}}+\log \log \left(\left(\frac{\sigma_q^{p+\nu}}{L^{q}}\right)^{\frac{1}{p+\nu-q}} \frac{1}{\varepsilon}\right)\right)\quad \text{if} \; \; q<p+\nu,
\end{gather*}
\begin{gather*}\label{complexity3}
\mathcal{O}\left(\left(\frac{L}{\sigma_q}\right)^{\frac{2}{3(p+\nu)-2}}\left(\frac{\sigma_q}{\varepsilon}\right)^{\frac{2(q-p-\nu)}{q(3(p+\nu)-2)}}\right) \quad \text{if} \;\; q>p+\nu.
\end{gather*}
\end{teo}

\subparagraph{The Proposed Approach:}
We are now in a position to combine all the building blocks and present our proposed approach. 
For solving the considered saddle point problem \eqref{main_spp_problem}, we reduce it to a combination of two auxiliary optimization problems separately for each group of variables: outer and inner problems, which are \eqref{equivalent_min_problem} and \eqref{inner_max_problem} respectively. By using Algorithms \ref{method_Pavel_2017} and \ref{restart_UAM} we propose the following approach.

\textbf{Approach 1:} 
The outer  minimization problem \eqref{equivalent_min_problem} is solved via ``Adaptive Gradient Method for Problems with Inexact Oracle'' (Algorithm \ref{method_Pavel_2017}), and in each iteration of Algorithm \ref{method_Pavel_2017}, the inner  maximization problem \eqref{inner_max_problem} is solved via ``Restarted Unified Acceleration Framework'' (Algorithm \ref{restart_UAM}).

Recall that the objective function $S(x,y) = F(x,y) - h(y) + r(x), \forall x\in \mathcal{X}, y \in \mathcal{Y} = \mathbb{R}^{n_y}$,  in the considered problem \eqref{main_spp_problem}, is non-convex in $x$, and uniformly concave in $y$ via the uniform convexity of $h$. The function $F$ has $(1, \nu)$-H\"older-continuous derivatives with respect to the primal variable $x$ and it has $(p, \nu)$-H\"older-continuous derivatives with respect to the dual variable $y$, for $\nu \in [0,1], p \in \{1,2, \ldots\}$, (see \eqref{smooth_F_1}---\eqref{smooth_F_4}). Let $\widehat{L} = \max\{L_{yx}, L_{yy}\}$.

For the complexity of the proposed approach, we find that from \eqref{complexity_alg_adaptive_grad}, the Algorithm \ref{method_Pavel_2017} will perform $\mathcal{O}\left( \frac{(L_{\nu}(g))^{\frac{1}{\nu_g}}\Delta}{\varepsilon^{\frac{1+\nu_g}{2 \nu_g}}}\right)$ steps (first-order oracle), where  $\Delta = g_0 - g_*$ is the difference between the value of the function $g$ at the initial point and its minimal value, and $L_{\nu}(g), \nu_g$ are given in Lemma \ref{lemma_function_g}. At each step of Algorithm \ref{method_Pavel_2017}, Algorithm \ref{restart_UAM} will perform at most the following number of iterations, 
dependently on $p, \nu$ and $q$ (the degree of the uniform convexity of the function $h$): 
\begin{gather*}
\mathcal{O}\left(\left(\frac{\widehat{L}}{\sigma_q}\right)^{\frac{2}{3(p+\nu)-2}} \log \left(\frac{1}{\varepsilon}\right) \right) \quad \text{if} \;\; q=p+\nu,
\end{gather*}
\begin{gather*}
\mathcal{O}\left(\left(\frac{\widehat{L}}{\sigma_q}\right)^{\frac{2}{3(p+\nu)-2}}+\log \log \left(\frac{1}{\varepsilon}\left(\frac{\sigma_q^{p+\nu}}{{\widehat{L}}^{q}}\right)^{\frac{1}{p+\nu-q}} \right)\right)\quad \text{if} \; \; q<p+\nu,
\end{gather*}
\begin{gather*}
\mathcal{O}\left(\left(\frac{\widehat{L}}{\sigma_q}\right)^{\frac{2}{3(p+\nu)-2}}\left(\frac{\sigma_q}{\varepsilon}\right)^{\frac{2(q-p-\nu)}{q(3(p+\nu)-2)}}\right) \quad \text{if} \;\; q>p+\nu.
\end{gather*}

From the previous,  we can see that the proposed approach is universal, in the sense of the order of smoothness of the objective function, where the considered class of problems contains more classes of problems as a special case, dependently on the values of the parameters $p, \nu$, and  $q$.

Summarizing, we have the following theorem, which gives the complexity of the proposed approach for solving the structured saddle point problem under consideration in the case when $q$, i.e., the degree of uniform concavity of the objective function \eqref{main_spp_problem}, equals $p + \nu$, which gives the better complexity.

\begin{teo}\label{complexity_approach}
By applying the proposed Approach 1, in the case when $q = p + \nu$,  we obtain an $\varepsilon$-solution of the problem \eqref{equivalent_min_problem} (which is equivalent to the problem \eqref{main_spp_problem}),  after $\mathcal{O}\left( \frac{(L_{\nu}(g))^{\frac{1}{\nu_g}}\Delta}{\varepsilon^{\frac{1+\nu_g}{2 \nu_g}}}\right)$ calls of the first-order oracle and 
$$
 \widetilde{\mathcal{O}}\left(\frac{(L_{\nu}(g))^{\frac{1}{\nu_g}}\Delta}{\varepsilon^{\frac{1+\nu_g}{2 \nu_g}}} \cdot \left(\frac{\widehat{L}}{\sigma_q}\right)^{\frac{2}{3(p+\nu)-2}} \right)
$$
calls of higher-order oracle, where $\widetilde{\mathcal{O}}(\cdot) = \mathcal{O}(\cdot)$ up to a logarithmic factor in $\log(\varepsilon^{-1})$.
\end{teo}

\paragraph{4. Conclusions}
In this paper, we developed a numerical method for solving a more general class of saddle point problems with composite structure and H\"older-continuous gradients, in the non-convex-uniformly-concave setting. We reduce the considered problem to a combination of two auxiliary optimization problems separately for each group of variables, namely, outer minimization problem (which is solved by the Adaptive Gradient Method applicable for non-convex problems and works with inexact oracle \cite{Dvurechensky_2017}) and inner maximization problem (which is solved by the Restarted Unified Acceleration Framework, which unifies high-order acceleration methods for minimizing a  convex function that has  H\"older-continuous higher-order derivative \cite{song2019unified}). We provided separate bounds for the number of calls to oracles for the outer and inner problems. Moreover, the complexity of the whole proposed approach is then estimated. As a future work it is planned to study the similar considered class of problems under the Polyak-Lojasiewicz condition instead of the non-convexity for the function with respect to the primal variable for the minimization problem, also to generalize the proposed approach for the non-Euclidean setting, and conduct numerical results for some applications and compare the results with other known approaches.

\end{document}